\documentclass[german,a4paper,12pt]{article}
\usepackage{amsmath,amsthm,amsfonts,latexsym,amscd}

\usepackage{amsthm}

\usepackage{amssymb}
\swapnumbers
\theoremstyle{plain}

\theoremstyle{definition}

\theoremstyle{remark}
\begin{document}

{\bf On the Newton filtration for functions on complete intersections}\\

{\bf Helmut A. Hamm (M\"unster}\\

A Newton filtration has been already considered by A.G. Kushnirenko [K] when he derived a formula for the Milnor number.\\
Let $g\in{\cal O}_{\mathbb{C}^n,0}$ be convenient (``commode''), i.e. $g$ contains for every $j\in\{1,\ldots,n\}$ at least one monomial of the form $z_j^l$. Assume $g(0)=0$.\\
Let $\Delta$ be the corresponding Newton polyhedron, i.e. the convex hull of $(supp\,g)+\mathbb{R}_+^n$. Let $p_1,\ldots,p_r\in \mathbb{N}^n$ be orthogonal to the compact facets (i.e. one-codimensional faces) of $\Delta$, the components being relatively prime. Let $\nu_{j}$ be maximal such that $<p_j,q>\ge \nu_j$ for all monomials $z^q$ which occur in $g$. This means that the facets are contained in the hypersurfaces $<p_j,q>= \nu_j$.\\
Choose $M\in\mathbb{N}\setminus\{0\}$ minimal such that $\psi:\mathbb{N}^n\to\mathbb{R}$: $\psi(q):=M\,\min\{\frac{1}{\nu_j}<p_j,q>\}$ takes values in $\mathbb{Z}$.\\
Then let $\hat{F}_l$ be the ideal in ${\cal O}_{\mathbb{C}^n,0}$ generated by all $z^q$ such that $\psi(q)\ge l$, i.e.$<p_j,q>\ge\frac{l}{M}\nu_j$ for all $j$.\\
This gives the classical (one-index) Newton filtration. Using this Kushnirenko proves a formula for the Milnor number in terms of the Newton polyhedron.\\

More recently there has been interest in multi-index filtrations. In particular, W.Ebeling and S.M.Gusein-Zade looked at a multi-index Newton filtration for functions on hypersurfaces [EG1]. As A.Lemahieu [L] showed, this may not be the induced filtration (which she called embedded filtration) which is simpler, we will primarily concentrate on the induced filtration.\\

The goal of the present paper is to pass from hypersurfaces to complete intersections. One motivation is that Kushnirenko looked in fact at special complete intersections given by the partial derivatives. At the end we will derive results which are similar, using other techniques.\\

{\bf 1. Multi-index Newton filtrations via local rings}\\

Let $g_1,\ldots,g_k\in{\cal O}_{\mathbb{C}^n,0}$ be germs of functions which define a complete intersection $(Y,0)$ in $(\mathbb{C}^n,0)$. Let $\Delta$ be the corresponding Newton polyhedron, i.e. the convex hull of $(supp\,g_1\cdot\ldots\cdot g_k)+\mathbb{R}_+^n$. We assume that $g_1,\ldots,g_k$ are convenient (``commode'') in the sense of Kushnirenko. Let $p_1,\ldots,p_r\in \mathbb{N}^n$ be orthogonal to the compact facets, the components being relatively prime. Let $\nu_{jl}$ be maximal such that $<p_l,q>\ge \nu_{jl}$ for all monomials $z^q$ which occur in $g_j$, $\nu_j:=(\nu_{j1},\ldots,\nu_{jr})$.\\
For $\mu\in\mathbb{Z}^r$ let $F_\mu{\cal O}_{\mathbb{C}^n,0}$ be the ideal generated by all $z^q$ such that $<p_l,q>\ge\mu_l$ for all $l$. This defines a multi-index filtration; we have\\
$F_{\mu'}\subset F_\mu$ for all $\mu\le\mu'$, i.e. with $\mu_l\le\mu'_l$ for all $l$,\\
$F_\mu\cdot F_{\mu'}\subset F_{\mu+\mu'}$.\\
Let ${\bf 1}:=(1,\ldots,1)$, then $\dim\,F_\mu/F_{\mu+{\bf 1}}<\infty$. Let us introduce the series\\
$L(t):= \dim\,F_\mu/F_{\mu+{\bf 1}} \,t^\mu$, where $t=(t_1,\ldots,t_r)$, and the Poincar\'e series\\
$P(t):=\frac{(t_1-1)\cdots(t_r-1)}{t_1\cdots t_r-1}L(t)$.\\
The reason is the following: if $V=\oplus W_\mu$ is a graded vector space and $F_\mu=\oplus W_{\mu'}$ where $\mu'\ge \mu$ we have $P(t)=\sum \dim W_\mu t^\mu$. See [CDG].\\
Therefore, the Poincar\'e series in the present case is $P_{{\cal O}_{\mathbb{C}^n,0}}(t)=\prod_{l=1}^n\frac{1}{1-t^{p^l}}$ where $p^l:=(p_{1l},\ldots,p_{rl})$.\\

Now the multi-index filtration on ${\cal O}_{\mathbb{C}^n,0}$ induces a multi-index filtration on ${\cal O}_{Y,0}$: $F_\mu{\cal O}_{Y,0}:=$ image of $F_\mu{\cal O}_{\mathbb{C}^n,0}$.\\
If $g_i=\sum c_{iq}z^q$ let $in_{j_1,\ldots,j_s}g_i$ be the initial part $\sum c_{iq}z^q$ where we take only those $q$ with $<p_{j_1},q>=\nu_{ij_1},\ldots, <p_{j_s},q>=\nu_{ij_s}$.\\

{\bf Theorem 1.1:} Assume that for all $s\ge 1,1\le j_1<\ldots<j_s\le r$ the following holds:\\
$in_{j_1,\ldots,j_s}g_1,\ldots,in_{j_1,\ldots,j_s}g_k$ define a subspace of $(\mathbb{C}^n,0)$ of codimension $\ge k-s+1$. Then\\
$P_{{\cal O}_{Y,0}}(t)=(1-t^{\nu_1})\cdots(1-t^{\nu_k}) P_{{\cal O}_{\mathbb{C}^n,0}}(t)$.\\

{\bf Special cases:} a) $k=1$: Only $s=1$ is relevant but the hypothesis is fulfilled since the initial part is non-trivial. So the only assumption is $g\neq 0$. See Lemahieu [L].\\
b) $k=2$: In the case $s=1$ we have a weak non-degeneracy condition. The case $s=2$, however, involves a heavy assumption: we need for all $j_1<j_2$ that
$in_{j_1,j_2}g_1$ or $in_{j_1,j_2}g_2$ is non-trivial. This implies that the corresponding facets of the Newton polyhedron of $g_1$ or those of $g_2$ intersect. If this holds for $g_1$ and all $j_1<j_2$ this means that the Newton polyhedron of $g_1$ is ``bistellar'', cf. Lemahieu.\\

The important step in the proof is the following:\\

{\bf Lemma 1.2:} Under the hypothesis of Theorem 1.1, the Koszul complex
$$0\to {\cal O}_{\mathbb{C}^n,0}\to \oplus {\cal O}_{\mathbb{C}^n,0}\to\ldots\to \oplus {\cal O}_{\mathbb{C}^n,0}\to{\cal O}_{\mathbb{C}^n,0}\to {\cal O}_{Y,0}\to 0$$
induces for every $\mu$ an exact sequence
$$0\to F_{\mu-\nu_1-\ldots-\nu_k}{\cal O}_{\mathbb{C}^n,0}\to \oplus F_{\mu-\nu_1-\ldots-\hat{\nu}_j-\ldots-\nu_k}{\cal O}_{\mathbb{C}^n,0}$$
$$\to\ldots\to \oplus F_{\mu-\nu_j}{\cal O}_{\mathbb{C}^n,0}\to F_\mu{\cal O}_{\mathbb{C}^n,0}\to F_\mu{\cal O}_{Y,0}\to 0.$$

{\bf Proof:} (i) First we deal with the exactness of\\
$0\to F_{\mu-\nu_1-\ldots-\nu_k}{\cal O}_{\mathbb{C}^n,0}\to \oplus F_{\mu-\nu_1-\ldots-\hat{\nu}_j-\ldots-\nu_k}{\cal O}_{\mathbb{C}^n,0}$\\
$\to\ldots\to \oplus F_{\mu-\nu_j}{\cal O}_{\mathbb{C}^n,0}\to F_\mu{\cal O}_{\mathbb{C}^n,0}$\hfill(*)\\
We pass from $\mathbb{C}^n$ to $\mathbb{C}^{n+r}$ with coordinates $z_1,\ldots,z_n,\zeta_1,\ldots,\zeta_r$. 
We have a $(\mathbb{C}^*)^r$-action on ${\cal O}:={\cal O}_{\mathbb{C}^{n+r},0}$:\\
$c\cdot z^q\zeta^\rho:=c_1^{<q,p_1>-\rho_1}\cdots c_r^{<q,p_r>-\rho_r}z^q\zeta^\rho$.\\
Let ${\cal O}(\mu)$ be the part where $c$ acts by multiplication by $c^\mu$: ``multiquasihomogeneous functions of (multi)degree $\mu$''. It consists of series involving monomials of the form\\
$z^q\zeta_1^{<p_1,q>-\mu_1}\cdots\zeta_r^{<p_r,q>-\mu_r}$. Note that $F_\mu{\cal O}_{\mathbb{C}^n,0}\simeq {\cal O}(\mu)$.\\
We replace $g_j\in F_{\nu_j}{\cal O}_{\mathbb{C}^n,0}$ by the corresponding element $G_j\in{\cal O}(\nu_j)$:
replace each monomial $z^q$ by\\
$z^q\zeta_1^{<p_1,q>-\nu_{j1}}\cdots\zeta_r^{<p_r,q>-\nu_{jk}}$. 

\vspace{2mm}
In order to prove (*) it is sufficient to show that $G_1,\ldots,G_k$ define a complete intersection in $(\mathbb{C}^{n+r},0)$:\hfill(**)\\
Then we obtain an exact Koszul complex
$$0\to {\cal O}\to \oplus  {\cal O}\to\ldots\to \oplus  {\cal O}\to  {\cal O}\to  {\cal O}/(G_1,\ldots,G_k)\to 0$$
This is an equivariant exact sequence: for each $\mu$, it gives rise to an exact sequence\\
$0\to {\cal O}(\mu-\nu_1-\ldots-\nu_k)\to \oplus  {\cal O}(\mu-\nu_1-\ldots-\hat{\nu}_j-\ldots-\nu_k)$\\
$\to\ldots\to \oplus {\cal O}(\mu-\nu_j)\to {\cal O}(\mu)\to  ({\cal O}/(G_1,\ldots,G_k))(\mu)\to 0$\\
Substituting $\zeta_1=\ldots=\zeta_r=1$ we obtain the desired exact sequence (*).\\
In fact, we have commutative diagrams of the form
$$\begin{array}{ccc}
{\cal O}(\mu-\nu_j)&\stackrel{G_j}{\longrightarrow}&{\cal O}(\mu)\\
\downarrow\simeq&&\downarrow\simeq\\
F_{\mu-\nu_j}{\cal O}_{\mathbb{C}^n,0}&\stackrel{g_j}{\longrightarrow}&F_\mu{\cal O}_{\mathbb{C}^n,0}
\end{array}$$
Now we prove (**). It is sufficient to show that $G_1=\ldots=G_k=0$ defines in $\{\zeta_{j_1}=\ldots=\zeta_{j_s}=0,\zeta_j\neq 0$ for $j\neq j_1,\ldots,j_s\}$ a subset of dimension $\le n+r-k$, i.e. of codimension $\ge k-s$. \\
The subset is invariant under the action of the torus $(\mathbb{C}^*)^r$, so we have to show the estimate for the codimension with $\zeta_j=1$ for $j\neq j_1,\ldots,j_s$. But this is true for $s=0$ because $g_1,\ldots,g_k$ define a complete intersection, for $s\ge 1$ by hypothesis of our Lemma (the codimension being even $\ge k-s+1$). So we have (**).

\vspace{2mm}
(ii) It remains to show the exactness of 
$$\oplus F_{\mu-\nu_j}{\cal O}_{\mathbb{C}^n,0}\to F_\mu{\cal O}_{\mathbb{C}^n,0}\to F_\mu{\cal O}_{Y,0}$$
Let $f\in F_\mu{\cal O}_{\mathbb{C}^n,0}$, $f\mapsto 0$ in $F_\mu{\cal O}_{Y,0}$. Then $f=\sum h_jg_j$. Choose $\mu'\le\mu$ maximal such that $h_l\in F_{\mu'-\nu_l}{\cal O}_{\mathbb{C}^n,0}$, $l=1,\ldots,k$. We have to show $\mu'=\mu$. Otherwise, $\mu'<\mu$, so there is a $j$ with $\mu'_j<\mu_j$. We have $f=\sum h_jg_j=0$ in $F_{\mu'}{\cal O}_{\mathbb{C}^n,0}/F_{\mu'+e_j}{\cal O}_{\mathbb{C}^n,0}$ because $f\in F_\mu{\cal O}_{\mathbb{C}^n,0}\subset F_{\mu'+e_j}{\cal O}_{\mathbb{C}^n,0}$.\\
It is sufficient to show that there is an exact sequence \\
$0\to F_{\mu'-\nu_1-\ldots-\nu_k}{\cal O}_{\mathbb{C}^n,0}/F_{\mu'-\nu_1-\ldots-\nu_k+e_j}{\cal O}_{\mathbb{C}^n,0} \to\ldots$\\
$\to F_{\mu'}{\cal O}_{\mathbb{C}^n,0}/F_{\mu'+e_j}{\cal O}_{\mathbb{C}^n,0}$ \hfill(***)\\
The last part of this sequence is 
$$\oplus_{l<m} F_{\mu'-\nu_l-\nu_m}{\cal O}_{\mathbb{C}^n,0}/F_{\mu'-\nu_l-\nu_m+e_j}{\cal O}_{\mathbb{C}^n,0}$$
$$\to \oplus_{l} F_{\mu'-\nu_l}{\cal O}_{\mathbb{C}^n,0}/F_{\mu'-\nu_l+e_j}{\cal O}_{\mathbb{C}^n,0}\to F_{\mu'}{\cal O}_{\mathbb{C}^n,0}/F_{\mu'+e_j}{\cal O}_{\mathbb{C}^n,0}$$
So in this case $h=(h_1,\ldots,h_k)\in \oplus F_{\mu'-\nu_l}{\cal O}_{\mathbb{C}^n,0}/F_{\mu'-\nu_l+e_j}{\cal O}_{\mathbb{C}^n,0}$ has an inverse image in
$\oplus_{l<m} F_{\mu'-\nu_l-\nu_m}{\cal O}_{\mathbb{C}^n,0}/F_{\mu'-\nu_l-\nu_m+e_j}{\cal O}_{\mathbb{C}^n,0}$ which can be used in order to replace $h$ by $\tilde{h}$ such that $\tilde{h}_l\in F_{\mu'+e_j-\nu_l}{\cal O}_{\mathbb{C}^n,0}$, in contradiction to the maximality of $\mu'$.\\
Now (***) is proved in an analogous way as (*) but here we really use the hypothesis of the Lemma.\\
More precisely: we have $F_{\mu'}{\cal O}_{\mathbb{C}^n,0}/F_{\mu'+e_j}{\cal O}_{\mathbb{C}^n,0}\simeq {\cal O}(\mu')/(\zeta_j)$. Note that we have commutative diagrams
$$\begin{array}{ccc}
{\cal O}(\mu'+e_j)&\stackrel{\cdot\zeta_j}{\longrightarrow}&{\cal O}(\mu')\\
\downarrow\simeq&&\downarrow\simeq\\
F_{\mu'+e_j}{\cal O}_{\mathbb{C}^n,0}&\longrightarrow&F_{\mu'}{\cal O}_{\mathbb{C}^n,0}
\end{array}$$
Of course, ${\cal O}\stackrel{\cdot\zeta_j}{\longrightarrow}{\cal O}$ is equivariant. 
Without loss of generality let $j=r$. We may replace $g_i$ by $in_rg_i$ in (***). Then we repeat the construction above with $r-1$ instead of $r$, and $p_1,\ldots,p_{r-1}$. We replace $in_rg_i$ by multiquasihomogeneous functions $\tilde{G}_i$, we must show that $\tilde{G}_1,\ldots,\tilde{G}_k$ define a complete intersection in $(\mathbb{C}^{n+r-1},0)$. As above we need that for $1<j_1<\ldots<j_s\le r-1$, $in_{j_1,\ldots,j_s,r}g_1=\ldots=in_{j_1,\ldots,j_s,r}g_k=0$ defines a subset of codimension $\ge k-s$. This is true by hypothesis.\\

{\bf Remark:} The structure of the proof may be more transparent if we first stick to a baby version of Lemma 1.2. Define\\
$\bar{F}_\mu:=F_\mu{\cal O}_{\mathbb{C}^n,0}/g_1F_{\mu-\nu_1}{\cal O}_{\mathbb{C}^n,0}+\ldots+g_kF_{\mu-\nu_k}{\cal O}_{\mathbb{C}^n,0}$.\\
The difficulty is that $\bar{F}_{\mu+{\bf 1}}\to\bar{F}_\mu$ might not be injective but it has a finite index. We can take the negative of the index in order to define the series $L$ and $P$ and have:\\

{\bf Lemma 1.3:} Assume that for all $1\le j_1<\ldots<j_s\le r$ the following holds:\\
$in_{j_1,\ldots,j_s}g_1,\ldots,in_{j_1,\ldots,j_s}g_k$ define a subspace of $(\mathbb{C}^n,0)$ of codimension $\ge k-s$.\\
Then we have an exact sequence
$$0\to F_{\mu-\nu_1-\ldots-\nu_k}{\cal O}_{\mathbb{C}^n,0}\to \oplus F_{\mu-\nu_1-\ldots-\hat{\nu}_j-\ldots-\nu_k}{\cal O}_{\mathbb{C}^n,0}\to\ldots$$
$$\to \oplus F_{\mu-\nu_j}{\cal O}_{\mathbb{C}^n,0}\to F_\mu{\cal O}_{\mathbb{C}^n,0}\to F_\mu{\cal O}_{\mathbb{C}^n,0}/(g_1F_{\mu-\nu_1}{\cal O}_{\mathbb{C}^n,0}+\ldots+g_kF_{\mu-\nu_k}{\cal O}_{\mathbb{C}^n,0})\to 0.$$
In particular, the Poincar\'e series for $\bar{F}$ is\\ 
$(1-t^{\nu_1})\cdots(1-t^{\nu_k}) P_{{\cal O}_{\mathbb{C}^n,0}}(t)$.\\

{\bf Proof:} Here we can content ourselves with the proof of (*) in Lemma 1.2.\\

{\bf Proof of Theorem 1.1:} We get a corresponding exact sequence as in Lemma 1.2 with $F_\mu$ replaced by $F_{\mu+{\bf 1}}$ or by $F_\mu/F_{\mu+{\bf 1}}$, similarly for the other indices. This yields the statement of Theorem 1.1 with $L$ instead of $P$. Now we can conclude for $P$.\\

The filtration which was considered in Theorem 1.1 was introduced by Lema\-hieu [L]. The filtration considered by Ebeling and Gusein-Zade [EG1] is different: it coincides with the one above if $\mu$ is a multiple of a standard unit vector $e_j$, in general $F_\mu$ is replaced by $F_{\mu_1e_1}\cap\ldots\cap F_{\mu_re_r}$. The question answered by Lemahieu in the hypüersurface case is under which condition both filtrations coincide.\\

{\bf Theorem 1.4:} Assume furthermore that for $s\ge 2, 0<j_1<\ldots<j_s\le r$ the following holds:\\
$in_{j_1,\ldots,j_s}g_1,\ldots,in_{j_1,\ldots,j_s}g_k$ define a subspace of codimension $\ge k-s+2$ in $(\mathbb{C}^n,0)$.\\
Then we have $F_\mu=F_{\mu_1e_1}\cap\ldots\cap F_{\mu_re_r}$ on ${\cal O}_{Y,0}$.\\

{\bf Special case:} $k=1$. Then only $s=2$ is relevant, and we need that $in_{j_1,j_2}g\neq 0$ for all $j_1<j_2$. This holds if and only if the Newton polyhedron is bistellar, as predicted by Lema\-hieu [L].\\

{\bf Proof of Theorem 1.4:} We may assume $\mu\in\mathbb{N}^r$. It is sufficient to look at the case where $\mu$ is not a scalar multiple of a standard unit vector, so $\mu$ is of the form $\mu'+e_j+e_l$ with $\mu'\in\mathbb{N}^r, j<l$. By induction it is sufficient to show that the induced filtration on ${\cal O}_{Y,0}$ satisfies the condition $F_{\mu+e_j+e_l}=F_{\mu+e_j}\cap F_{\mu+e_l}$.\\
This condition is equivalent to the exactness of the sequence
$$0\to F_{\mu+e_j+e_l}\stackrel{h\mapsto (h,h)}{\longrightarrow} F_{\mu+e_j}\oplus F_{\mu+e_l} \stackrel{(h_1,h_2)\mapsto h_1-h_2}{\longrightarrow} F_{\mu+e_j}+F_{\mu+e_l}\to 0$$
So we have to show that the Koszul resolution induces an exact sequence for $F_{\mu'+e_j}+F_{\mu'+e_l}$ or for $F_{\mu'}/(F_{\mu'+e_j}+F_{\mu'+e_l})$ (with suitable choices of $\mu'$). This is shown in a similar way as (*) or (***) in the proof of Lemma 1.\\
In fact, $F_{\mu'}/(F_{\mu'+e_j}+F_{\mu'+e_l})\simeq {\cal O}(\mu')/(\zeta_j,\zeta_l)$. Note that we have commutative diagrams with exact rows:
\footnotesize
$$\begin{array}{ccccccc}
{\cal O}(\mu'+e_j)\oplus{\cal O}(\mu'+e_l)&\stackrel{(\cdot\zeta_j,\cdot\zeta_l)}{\longrightarrow}&{\cal O}(\mu')&\to&
{\cal O}(\mu')/(\zeta_j,\zeta_l)&\to&0\\
\downarrow\simeq&&\downarrow\simeq&&\downarrow&&\\
F_{\mu+e_j}{\cal O}_{\mathbb{C}^n,0}\oplus F_{\mu+e_l}{\cal O}_{\mathbb{C}^n,0}&\longrightarrow&F_{\mu'}{\cal O}_{\mathbb{C}^n,0}&\to&F_{\mu'}{\cal O}_{\mathbb{C}^n,0}/(F_{\mu'+e_j}{\cal O}_{\mathbb{C}^n,0}+F_{\mu'+e_l}{\cal O}_{\mathbb{C}^n,0})&\to&0
\end{array}$$
\normalsize
Suppose without loss of generality that $j=r-1,l=r$. Then we work with $\mathbb{C}^{n+r-2}$ instead of $\mathbb{C}^{n+r}$. Replace $in_{r-1,r}g_i$ by a multiquasihomogeneous $\hat{G}_i$. We must show that $\hat{G}_1=\ldots=\hat{G}_k=0$ defines a complete intersection in $(\mathbb{C}^{n+r-2},0)$, which follows from the hypothesis (with $j_{s-1}=r-1,j_s=r$).\\ 

Now we compare with the result of Kushnirenko [K]. There $k=n,\nu_1=\ldots=\nu_n=\nu$, $\mu=\frac{l}{d}\nu$ with $l\in\mathbb{N}$. If we specialize Theorem 1.1 to this case we get a much weaker result than Kushnirenko since he does not impose any condition on the Newton polyhedron, the essential point is a non-degeneracy condition.\\

The disadvantage in the approach above is that from the Newton polyhedron we used only a partial information, concerning the facets. So we will take toric varieties into account.\\

{\bf 2. Multi-index Newton filtrations via toric varieties}\\

Let $g_1,\ldots,g_k,\Delta$ be as in \S 1. We may use a toric variety $X$ associated with the Newton polyhedron $\Delta$, cf. [KKMS] or [O]. Let us start with the corresponding fan: it is a collection of cones whose union is $\mathbb{R}_+^n$. The edges (i.e. one-dimensional faces) are generated by $e_1,\ldots,e_n,p_1,\ldots,p_r$. The additional information is which of them span a cone of the fan. Choosing some simplicial subdivision we replace it by a new fan. We have a corresponding toric variety $X$ and a proper birational mapping $\pi:X\to\mathbb{C}^n$. Let $D:=\pi^{-1}(\{0\})$ be the exceptional divisor.\\

From now on we work in the complex analytic category, otherwise we write $X^{alg}$ instead of $X$.\\

Now $X^{alg}$ is the union of affine toric varieties $X^{alg}_\sigma$, where $\sigma$ is a cone of the fan. Note that ${\cal O}(X^{alg}_\sigma)$ is spanned by all $z^q$, $q\in\mathbb{Z}^n$, such that $q_{i_1}\ge 0,\ldots,q_{i_l}\ge 0, <p_{j_1},q>\ge 0,\ldots,<p_{j_s},q>\ge 0$ if $\sigma$ is spanned by $e_{i_1},\ldots,e_{i_l},p_{j_1},\ldots,p_{j_s}$. Let ${\cal O}_{X^{alg}}<D>:=j_*{\cal O}_{X^{alg}\setminus D^{alg}}$, where $j:X^{alg}\setminus D^{alg}\to X^{alg}$ is the inclusion. If $\sigma$ is as above, we have that  $H^0(X^{alg}_\sigma,{\cal O}_{X^{alg}}<D>)$ is generated by all $z^q$, $q\in\mathbb{Z}^n$, such that $q_{i_1},\ldots,q_{i_l}\ge 0$. Let ${\cal O}_X<D>:=({\cal O}_{X^{alg}}<D>)^{an}$ be the sheaf of meromorphic functions on $X$ which are holomorphic on $X\setminus D$. For $\mu\in\mathbb{Z}^r$ we can introduce ${\cal O}_X$-submodules ${\cal F}_\mu$ of ${\cal O}_X<D>$: if $\sigma$ is as above, let ${\cal F}_\mu(X_\sigma)$ be generated by all $z^q$, $q\in\mathbb{Z}^n$, such that $q_{i_1}\ge 0,\ldots,q_{i_l}\ge 0, <p_{j_1},q>\ge \mu_{j_1},\ldots,<p_{j_s},q>\ge \mu_{j_s}$. So we have a filtration ${\cal F}_\mu$ on ${\cal O}_X<D>(X_\sigma)$. Note that ${\cal O}_X<D>$ is the union of all ${\cal F}_\mu$, and the sheaves ${\cal F}_\mu$ are equivariant with respect to the action of the torus $(\mathbb{C}^*)^n$.\\
Define $h_\mu:\mathbb{R}_+^n\to\mathbb{R}$ as follows: $h_\mu(e_i):=0$ for $i=1,\ldots,n$, $h_\mu(p_j):=\mu_j$ for $j=1,\ldots,r$, $h_\mu|\sigma$ linear for all cones $\sigma$ of the fan. Then $h_\mu$ is well-defined because the fan was assumed to be simplicial. Now ${\cal F}_\mu(X_\sigma)$ is generated by all $z^q$, $q\in\mathbb{Z}^n$, such that $<p,q>\ge h_\mu(p)$ for all $p\in\sigma$. Note that $h_\mu$ does not have to be a support function in the sense of [O] p.66 since we cannot guarantee that $h_\mu(\mathbb{N}^n)\subset\mathbb{Z}$. In particular, ${\cal F}_\mu$ does not need to be invertible. Note, however, that for $\mu=\nu_{l_1}+\ldots+\nu_{l_s}$, $1\le l_1<\ldots<l_s\le k$, $h_\mu$ is convex, $h_\mu(\mathbb{N}^n)\subset \mathbb{Z}$, so ${\cal F}_\mu$ is invertible.\\ 

Let $\hat{g}_j:=g_j\circ \pi$, it is defined on a neighbourhood of $D$ in $X$.\\

We want to restrict to conditions which hold for a generic choice of the coefficients of the $g_i$.\\

For every cone $\sigma$ of the fan, spanned by $e_{i_1},\ldots,e_{i_l},p_{j_1},\ldots,p_{j_s}$, let $in_\sigma g_i$ be the part involving only those $z^q$ such that $q_{i_1}=\ldots=q_{i_l}=0, <p_{j_1},q>=\nu_{ij_1},\ldots, <p_{j_s},q>=\nu_{ij_s}$.\\

Let $\Delta_i$ be the convex hull of $(supp\,g_i)+\mathbb{N}^n$, $\Gamma_i$ the union of the compact faces of $\Delta_i$, $S_i:=\Gamma_i\cap\mathbb{N}^n$, $i=1,\ldots,k$. Note that $S_i$ is finite.\\

{\bf Lemma 2.1:} Suppose that the coefficients of the functions $g_i$ at $z^q, q\in S_i$, are chosen generically. Then $g_1=\ldots=g_k=0$ defines in $(\mathbb{C}^n,0)$ a complete intersection with isolated singularity, and for every cone $\sigma$ of the fan, spanned by $e_{i_1},\ldots,e_{i_l},p_{j_1},\ldots,p_{j_s}$, with $s\ge 1$ the functions $in_\sigma g_i$, $i=1,\ldots,k$, define a submanifold of $(\mathbb{C}^*)^n$ of codimension $k$.\\

{\bf Proof:} Define $G_i:\mathbb{C}^{S_1}\times\ldots\times\mathbb{C}^{S_k}\times(\mathbb{C}^*)^n\to\mathbb{C}$: $((a_{1q}),\ldots,(a_{kq}),z)\mapsto \sum_{q\in S_i}\, a_{iq}z^q$. We can define $in_\sigma G_i$ similarly as in the case of $g_i$. For each $\sigma$ as above,  $in_\sigma G_1=\ldots=in_\sigma G_k=0$ defines a submanifold of $\mathbb{C}^{S_1}\times\ldots\times\mathbb{C}^{S_k}\times(\mathbb{C}^*)^n$ of codimension $k$. We restrict the canonical projection onto $\mathbb{C}^{S_1}\times\ldots\times\mathbb{C}^{S_k}$ to this submanifold; the regular values for all $\sigma$ altogether yield the desired generic choice of coefficients.\\
Now fix $\sigma$. On $X_\sigma$ we may trivialize ${\cal F}_{\nu_i}$. We regard $\hat{g}_1,\ldots,\hat{g}_k$ as sections of ${\cal F}_\mu|X_\sigma$; taking a trivialization we can consider them as functions $\hat{\rm g}_1,\ldots,\hat{\rm g}_k$ on $X_\sigma$. The restriction of the latter to the minimal orbit defined by $\sigma$ is smooth of codimension $k$ by assumption. We get smoothness of $\{\hat{\rm g}_1=\ldots=\hat{\rm g}_k=0\}$ in a neighbourhood outside $D$ because the subdivision of $X$ into orbits is a Whitney-regular stratification. Altogether we obtain that $X\cap\{g_1=\ldots=g_k=0\}$ has an isolated singularity at $0$.\\

In the following lemma we can restrict ourselves to a weaker hypothesis.\\

{\bf Lemma 2.2:} Suppose that we have the following condition:\\
For every cone $\sigma$ of the fan, spanned by $e_{i_1},\ldots,e_{i_l},p_{j_1},\ldots,p_{j_s}$, $s\ge 1$, the functions $in_\sigma g_i$, $i=1,\ldots,k$, define a subspace of $(\mathbb{C}^*)^n$ of codimension $\ge k-l-s$.\\
Then we have an exact sequence:\\
$0\to {\cal F}_{\mu-\nu_1-\ldots-\nu_k}|D\to \oplus {\cal F}_{\mu-\nu_1-\ldots-\hat{\nu}_j-\ldots-\nu_k}|D\to\ldots\to\oplus {\cal F}_{\mu-\nu_j}|D\to {\cal F}_\mu|D\to ({\cal F}_{\mu}|D)/(\hat{g}_1{\cal F}_{\mu-\nu_1}|D+\ldots+\hat{g}_k{\cal F}_{\mu-\nu_k}|D)\to 0$.\\

{\bf Proof:} It is sufficient to prove the exactness for the restriction to $X_\sigma\cap D$, $\sigma$ as above (with $s\ge 1$, otherwise $X_\sigma\cap D=\emptyset$).\\
Without loss of generality assume that $(i_1,\ldots,i_l)=(1,\ldots,l)$ and $(j_1,\ldots,j_s)=(1,\ldots,s)$.\\
Let us look at the filtration $F_\mu^\sigma:=H^0(X_\sigma\cap D,{\cal F}_\mu)$ on $H^0(X_\sigma\cap D,{\cal O}_X<D>)$.\\
Now $0\to F_{\mu-\nu_1-\ldots-\nu_k}^\sigma\to \oplus F_{\mu-\nu_1-\ldots-\hat{\nu}_j-\ldots-\nu_k}^\sigma\to\ldots\to \oplus F_{\mu-\nu_j}^\sigma\to F_\mu^\sigma$\\
$\to F_\mu^\sigma/(g_1F_{\mu-\nu_1}^\sigma+\ldots+g_kF_{\mu-\nu_k}^\sigma)\to 0$ is exact:\\
Similarly as in the proof of Lemma 1.2 we introduce on \\
$A:={\cal O}_{\mathbb{C}^{n+r},0}[z_{l+1}^{-1},\ldots,z_n^{-1},\zeta_{s+1}^{-1},\ldots,\zeta_r^{-1}]$ the analogous $(\mathbb{C}^*)^r$-action. Let $A(\mu)$ be the subspace where $c$ acts as multiplication by $c^\mu$. Then we have: $F_\mu^\sigma\simeq A(\mu)$. Let $G_i$ be defined as in the proof of Lemma 1.2. Then we have to show that we have an exact sequence:\\
$0\to A(\mu-\nu_1-\ldots-\nu_k)\to\ldots \to\oplus A(\mu-\nu_j)\to A(\mu)\to (A/(G_1,\ldots,G_k))(\mu)\to 0$.\\
This follows from the exactness of the sequence\\
$0\to A\to\ldots \to\oplus A\to A\to A/(G_1,\ldots,G_k)\to 0$\\
which holds as soon as we show that $G_1,\ldots,G_k$ is a regular sequence on $A$.\\
It is sufficient to look at the orbits with respect to the torus action.\\
The minimal orbit is given by $z_1=\ldots=z_l=0, z_{l+1}\cdot\ldots\cdot z_n\neq 0$, $\zeta_1=\ldots=\zeta_s=0$, $\zeta_{s+1}\cdot\ldots\cdot \zeta_r\neq 0$. We have to show that $G_1=\ldots=G_k=0$ defines a subset of this orbit of dimension $\le n+r-k$, i.e. of codimension $\ge k-l-s$. We may pass to the subset given by $z_{l+1}=\ldots=z_n=\zeta_{s+1}=\ldots=\zeta_r=1$. Now we have to prove that $in_\sigma g_i$, $i=1,\ldots,k$, define a subset of the same codimension in $(\mathbb{C}^*)^{n-l}$ or also in $(\mathbb{C}^*)^n$. But this follows from the assumption.\\
For the other orbits it means to pass to faces of $\sigma$ because the fan of $X$ is supposed to be simplicial.\\

{\bf Lemma 2.3:} Suppose that we have the condition:\\
For every cone $\sigma$ of the fan, spanned by $e_{i_1},\ldots,e_{i_l},p_{j_1},\ldots,p_{j_s}$, $s\ge 1$, the functions $in_\sigma g_i$, $i=1,\ldots,k$, define a subspace of $(\mathbb{C}^*)^n$ of codimension $\ge k-l-s+1$. Then:\\
a)  $({\cal F}_{\mu+{\bf 1}}|D)/(\hat{g}_1{\cal F}_{\mu-\nu_1+{\bf 1}}+\ldots+\hat{g}_k{\cal F}_{\mu-\nu_k+{\bf 1}}|D)\to ({\cal F}_{\mu}|D)/(\hat{g}_1{\cal F}_{\mu-\nu_1}+\ldots+\hat{g}_k{\cal F}_{\mu-\nu_k}|D)$ is injective.\\
b) Let $(\hat{g}_1,\ldots,\hat{g}_k)$ be the ideal in ${\cal O}_X<D>$ generated by $\hat{g}_1,\ldots,\hat{g}_k$. Then we have:\\ 
$({\cal F}_{\mu}|D)/(\hat{g}_1{\cal F}_{\mu-\nu_1}|D+\ldots+\hat{g}_k{\cal F}_{\mu-\nu_k}|D)\simeq ({\cal F}_\mu|D)/((\hat{g}_1,\ldots,\hat{g}_k)\cap({\cal F}_\mu|D))$.\\

{\bf Proof:} a) With the notations before we have the following commutative diagram:
$$\begin{array}{ccc}
F_{\mu+{\bf 1}}^\sigma&\longrightarrow& F_\mu^\sigma\\
\downarrow\simeq&&\downarrow\simeq\\
A(\mu+{\bf 1})&\stackrel{\cdot \zeta_1\cdots\zeta_r}{\longrightarrow}&A(\mu)\\
\end{array}$$
So it is sufficient to show that the multiplication by $\zeta_1\cdot\ldots\cdot\zeta_r$ defines an injective mapping
$A/(G_1,\ldots,G_k)\to A/(G_1,\ldots,G_k)$,\\
or that $\zeta_1\cdot\ldots\cdot\zeta_r,G_1,\ldots,G_k$ define a regular sequence on $A$.\\
Here we have to look at the dimensions of sets and just need the hypothesis.\\
Alternative: use b).\\
b) Here we argue as in part (ii) of the proof of Lemma 1.2.\\

As a corollary we obtain:\\

{\bf Theorem 2.4:} Suppose that we have the hypothesis of Lemma 2.3.\\
a) We have an exact sequence:\\
$0\to {\cal F}_{\mu-\nu_1-\ldots-\nu_k}/{\cal F}_{\mu-\nu_1-\ldots-\nu_k+{\bf 1}}\to \oplus {\cal F}_{\mu-\nu_1-\ldots-\hat{\nu}_j-\ldots-\nu_k}/{\cal F}_{\mu-\nu_1-\ldots-\hat{\nu}_j-\ldots-\nu_k+{\bf 1}}\to\ldots\to \oplus {\cal F}_{\mu-\nu_j}/{\cal F}_{\mu-\nu_j+{\bf 1}}\to {\cal F}_\mu/{\cal F}_{\mu+{\bf 1}}\to {\cal F}_{\mu}/({\cal F}_{\mu+{\bf 1}}+\hat{g}_1{\cal F}_{\mu-\nu_1}+\ldots+\hat{g}_k{\cal F}_{\mu-\nu_k})\to 0$\\
b) $\tilde{P}(t)=(1-t^{\nu_1})\cdots(1-t^{\nu_k})P(t)$.\\

As for b), let $L(t):=\sum \chi(D,{\cal F}_{\mu}/{\cal F}_{\mu+{\bf 1}})t^\mu$, $\tilde{L}(t):= \sum \chi(D,{\cal F}_{\mu}/({\cal F}_{\mu+{\bf 1}}+\hat{g}_1{\cal F}_{\mu-\nu_1}+\ldots+\hat{g}_k{\cal F}_{\mu-\nu_k}))t^\mu$, and let $P, \tilde{P}$ be defined analogously. Of course, $\chi(D,{\cal S}):=\sum_i(-1)^i\dim H^i(D,{\cal S})$ if $\cal S$ is a coherent analytic sheaf on $D$.\\

{\bf Proof:} a) Because of Lemma 2.3a) we have a short exact sequence of complexes on $D$ (we omit to write the restriction to $D$):
\footnotesize
$$\begin{array}{ccccccc}
&0&&0&&0&\\
&\downarrow&&\downarrow&&\downarrow&\\
0\to&{\cal F}_{\mu-\nu_1-\ldots-\nu_k+{\bf 1}}&\to\ldots\to&{\cal F}_{\mu+{\bf 1}}&\to&{\cal F}_{\mu+{\bf 1}}/(\hat{g}_1{\cal F}_{\mu-\nu_1+{\bf 1}}+\ldots+\hat{g}{\cal F}_{\mu-\nu_k+{\bf 1}})&\to0\\
&\downarrow&&\downarrow&&\downarrow&\\
0\to&{\cal F}_{\mu-\nu_1-\ldots-\nu_k}&\to\ldots\to&{\cal F}_\mu&\to&{\cal F}_\mu/(\hat{g}_1{\cal F}_{\mu-\nu_1}+\ldots+\hat{g}{\cal F}_{\mu-\nu_k})&\to0\\
&\downarrow&&\downarrow&&\downarrow&\\
0\to&{\cal F}_{\mu-\nu_1-\ldots-\nu_k}/{\cal F}_{\mu-\nu_1-\ldots-\nu_k+{\bf 1}}&\to\ldots\to&{\cal F}_\mu/{\cal F}_{\mu+{\bf 1}}&\to&{\cal F}_\mu/({\cal F}_{\mu+{\bf 1}}+\hat{g}_1{\cal F}_{\mu-\nu_1}+\ldots+\hat{g}{\cal F}_{\mu-\nu_k})&\to0\\
&\downarrow&&\downarrow&&\downarrow&\\
&0&&0&&0&
\end{array}$$
\normalsize
Since the two first complexes are exact by Lemma 2.2 the last one is exact, too.\\
b) follows from a).\\

In order to apply Theorem 2.4 it is useful to develop a formula for $L(t)$. We introduce numbers $\chi_I$ and $n_{I,\mu}$:\\
Let $\sigma_1,\ldots,\sigma_l$ be the maximal cones of the fan which defines $X$. Put $I_{\mu,q}:=\{(i,j)\,|\,1\le i\le n, 1\le j\le r, q_i\ge 0,<p_j,q>\ge\mu_j\}$. For $\lambda=1,\ldots,l$ let $\sigma_\lambda$ be generated by the $e_i$, $i\in J_\lambda^1$, and $p_j$,  $j\in J_\lambda^2$. Put $J_\lambda:=J_\lambda^1\times J_\lambda^2$. Let $I=I^1\times I^2$ be given, $I^1\subset
\{1,\ldots,n\}, I^2\subset\{1,\ldots,r\}$. Let $m_{I,\nu}$ be the number of all $\Lambda\subset\{1,\ldots,l\}$ such that $\#\Lambda=\nu$ and $\bigcap_{\lambda\in \Lambda}J_\lambda\subset I$. Put $\chi_I:=\sum_{\nu>0} (-1)^{\nu-1}m_{I,\nu}$.\\
On the other hand, let $n_{I,J,\mu}$ be the number of all $q$ such that $I=I_{\mu,q}$ and $J= I_{\mu+{\bf 1},q}$.\\

{\bf Theorem 2.5:} $L(t)=\sum_\mu\sum_{J\subset I}\,n_{I,J,\mu}(\chi_I-\chi_J)t^\mu$.\\

We will see that $n_{I,J,\mu}<\infty$ as soon as $\chi_I\neq \chi_J$, so the sum is well-defined.\\

{\bf Proof:} For a cone $\sigma$ of the fan, generated by $e_{i_1},\ldots,e_{i_l}, p_{j_1},\ldots,p_{j_s}$, $X_\sigma$ is affine, so $\chi(X_\sigma,{\cal F}_\mu)(q)=\dim\,H^0(X_\sigma,{\cal F}_\mu)(q)=1$ if $<p_{j_1},q>\ge\mu_{j_1},\ldots,<p_{j_s},q>\ge\mu_{j_s},q_{i_1}\ge 0,\ldots,q_{i_l}\ge 0$; otherwise we get $0$. Since $X_{\sigma_1},\ldots,X_{\sigma_l}$ form an affine covering of $X$ and $X_{\sigma_{\lambda_1}}\cap\ldots\cap X_{\sigma_{\lambda_\nu}}= X_{\sigma_{\lambda_1}\cap\ldots\cap\sigma_{\lambda_\nu}}$ we have:\\ 
$\chi(X,{\cal F}_\mu)(q)=\sum_{\nu>0} (-1)^{\nu-1}\sum _{1\le\lambda_1<\ldots<\lambda_\nu\le l} \chi(X_{\sigma_{\lambda_1}\cap\ldots\cap\sigma_{\lambda_\nu}},{\cal F}_\mu)(q)=\\
\sum(-1)^{\nu-1}\#\{\Lambda\subset\{1,\ldots,l\}\,|\,\#\Lambda=\nu,\cap_{\lambda\in\Lambda}J_\lambda\subset I_{\mu,q}\}=\chi_{I_{\mu,q}}$.\\
So $\chi(X,{\cal F}_\mu/{\cal F}_{\mu+{\bf 1}})(q)=\chi_{I_{\mu,q}}-\chi_{I_{\mu+{\bf 1},q}}$.\\
Since ${\cal F}_\mu/{\cal F}_{\mu+{\bf 1}}$ is coherent on the compact space $D$ there are for each $\mu$ only finitely many $q$ such that $\chi_{I_{\mu,q}}\neq \chi_{I_{\mu+{\bf 1},q}}$. This implies our assertion.\\

Now let us return to our original question. Here we compare the filtration $F_\mu$ of \S 1 with the filtration given by $H^0(D,{\cal F}_\mu/(\hat{g}_1,\ldots,\hat{g}_k)\cap{\cal F}_\mu)$ on $H^0(D,{\cal O}_X<D>/(\hat{g}_1,\ldots,\hat{g}_k))$. So have to concentrate on $H^0$ instead of the Euler characteristic.\\

Let $Y'$ be the strict transform of $Y$ with respect to $\pi$ and $\pi':=\pi|Y'$, $D':=D\cap Y'$. Then $H^0(D',{\cal O}_{Y'}<D'>)\simeq H^0(D,{\cal O}_X<D>/(\hat{g}_1,\ldots,\hat{g}_k))$.\\

Recall that we are working in the complex analytic category: if we pass to the complex algebraic one we put an index ``alg''.\\

Let $D_i$ be the component of $D$ which corresponds to $\zeta_i=0$, $D'_i:=D_i\cap Y'$. If $D'_1,\ldots,D'_r$ are non-empty and irreducible we have a divisorial filtration on ${\cal O}_{Y,0}$ and ${\cal M}_{Y,0}$ with respect to $D'_1,\ldots,D'_r$, see [EG2], where ${\cal M}_Y$ denotes the sheaf of meromorphic functions: if $h\in{\cal O}_{Y,0}$ resp ${\cal M}_{Y,0}$, look at the order of $h\circ\pi'$ with respect to $D'_i$, $i=1,\ldots,r$. If we drop the assumption on non-emptiness and irreducibility, we look at the order of $h\circ\pi'$ with respect to all irreducible components of $D'_i$ instead and speak of the ``generalized divisorial filtration''.\\

Then we have:\\

{\bf Lemma 2.6:} Suppose that the following holds:\\
For every cone $\sigma$ of the fan, spanned by $e_{i_1},\ldots,e_{i_l},p_{j_1},\ldots,p_{j_s}$, $s\ge 1$, the functions $in_\sigma g_i$, $i=1,\ldots,k$, define a smooth complete intersection in $(\mathbb{C}^*)^n$.\\
a) Suppose $n-k\ge 2$. Then the spaces $H^0(D,{\cal F}_\mu/(\hat{g}_1,\ldots,\hat{g}_k)\cap {\cal F}_\mu))$ correspond to the generalized divisorial filtration on ${\cal O}_{Y,0}\simeq H^0(D,{\cal O}_X/(\hat{g}_1,\ldots,\hat{g}_k)\cap{\cal O}_X)=H^0(D,{\cal O}_X<D>/(\hat{g}_1,\ldots,\hat{g}_k))$ with respect to the $D_i'$.\\
b) Suppose $n-k=1$. Then $\pi':Y'\to Y$ is the normalisation of $Y$, and the spaces $H^0(D,{\cal F}_\mu/(\hat{g}_1,\ldots,\hat{g}_k)\cap {\cal F}_\mu))$ correspond to the generalized divisorial filtration on ${\cal M}_{Y,0}\simeq H^0(D',{\cal O}_{Y'}<D'>)$with respect to the $D_i'$\\
c) Suppose that $n-k=0$. Then $H^0(D,{\cal F}_\mu/(\hat{g}_1,\ldots,\hat{g}_k)\cap{\cal F}_\mu)=0$ for all $\mu$.\\

{\bf Proof:} As we have seen in Lemma 2.1, $(Y,0)$ is a complete intersection with isolated singularity.\\
a) Note that $Y$ is a complete intersection of dimension $\ge 2$ with an isolated singularity, hence normal by the Serre criterion for normality [H] II Theorem 8.22A.\\
By the proof of Zariski's Main Theorem [H] III Cor. 11.4, we have ${\cal O}_{Y,0}\simeq H^0(D',{\cal O}_{Y'})\simeq H^0(D,{\cal O}_X/(\hat{g}_1,\ldots,\hat{g}_k)\cap {\cal O}_X)$. Let $j:Y^{alg}\setminus\{0\}\to Y^{alg}$ be the inclusion. By normality, ${\cal O}_{Y,0}\simeq (j_*{\cal O}_{Y^{alg}\setminus\{0\}})_0^{an}\simeq H^0(D',{\cal O}_{Y'}<D'>)\simeq H^0(D,{\cal O}_X<D>/(\hat{g}_1,\ldots,\hat{g}_k))$.\\
Obviously the generalized divisorial filtration is as described.\\
b) Note that $Y'$ is smooth. Furthermore ${\cal M}_{Y,0}\simeq (j_*{\cal O}_{Y^{alg}\setminus\{0\}})_0^{an}\simeq H^0(D',{\cal O}_{Y'}<D'>)\simeq H^0(D,{\cal O}_X<D>/(\hat{g}_1,\ldots,\hat{g}_k))$.\\
c) We have that $D\cap\{\hat{g}_1=\ldots=\hat{g}_n=0\})=\emptyset$.\\

{\bf Remark:} a) Assume $k=1, n\ge 3$. Then the spaces $D'_i$ are irreducinble, see [EG2], so we can speak of a divisorial filtration. For a proof and further results, see the Appendix (section 4).\\
b) Assume furthermore $n-k=1$. Then the Poincar\'e series for the filtration given by $H^0(D,{\cal F}_\mu/(\hat{g}_1,\ldots,\hat{g}_k)\cap{\cal F}_\mu)$ is $\tilde{P}(t)$, see Theorem 2.4: Note that the sheaves in question are concentrated upon $D'$, and $\dim\,D'=0$, so $\dim\,H^0$ coincides with the Euler characteristic.\\

Now let us return to the filtration of \S 1.\\

{\bf Theorem 2.7:} Suppose that we have the assumption of Lemma 2.3 and that $h_{\mu-\nu_1-\ldots-\nu_k}$ is convex.\\
a) The sequence\\
$0\to F_{\mu-\nu_1-\ldots-\nu_k}/F_{\mu-\nu_1-\ldots-\nu_k+{\bf 1}}\to\ldots \to F_\mu/F_{\mu+{\bf 1}}\to \bar{F}_\mu\to 0$\\
with $\bar{F_\mu}:= F_\mu/(F_{\mu+{\bf 1}}+g_1F_{\mu-\nu_1}+\ldots+g_kF_{\mu-\nu_k})$ is exact.\\
b) The coefficients of the following series at $t^\mu$ coincide:\\
(i) $(1-t^{\nu_1})\cdots(1-t^{\nu_k})P_{{\cal O}_{\mathbb{C}^n,0}}(t)$,\\
(ii) Poincar\'e series for $\bar{F}_{\mu'}$,\\
(iii) Poincar\'e series for $H^0(D,{\cal F}_{\mu'}/(\hat{g}_1,\ldots,\hat{g}_k)\cap{\cal F}_{\mu'})$.\\

{\bf Remark:} a) Note that the relation between (i) and (ii) is not obvious because we cannot apply Lemma 1.3.\\
If $n>k$ the filtration in (iii) is the generalized divisorial filtration, see Lemma 2.6.\\

{\bf Proof:} a) First we show that there is an exact sequence:
$0\to F_{\mu-\nu_1-\ldots-\nu_k}\to\ldots\to F_\mu\to F_\mu/(g_1F_{\mu-\nu_1}+\ldots+g_kF_{\mu-\nu_k})\to 0$.\hfill(*)\\ 
Then we prove that we have the same exact sequence with $\mu+{\bf 1}$ instead of $\mu$. \\
Moreover we show the injectivity of \\
$F_{\mu+{\bf 1}}/(g_1F_{\mu-\nu_1+{\bf 1}}+\ldots+g_kF_{\mu-\nu_k+{\bf 1}})\to F_\mu/(g_1F_{\mu-\nu_1}+\ldots+g_kF_{\mu-\nu_k})$ \hfill(**).\\
Using a short exact sequence of sheaf complexes we obtain a).\\
It remains to prove (*), its analogue with $\mu+{\bf 1}$ instead of $\mu$, and (**).

\vspace{2mm}
In order to show (*) we prove that there is a commutative diagram with exact rows:

\footnotesize
$$\begin{array}{ccccccc}
0\to &F_{\mu-\nu_1-\ldots-\nu_k}&\to\ldots\to&F_\mu&\to&F_\mu/(g_1F_{\mu-\nu_1}+\ldots+g_kF_{\mu-\nu_k})&\to 0\\
&\downarrow\simeq&&\downarrow\simeq&&\downarrow\simeq&\\
0\to &H^0(D,{\cal F}_{\mu-\nu_1-\ldots-\nu_k})&\to\ldots\to&H^0(D,{\cal F}_\mu)&\to&H^0(D,{\cal F}_\mu/\hat{g}_1{\cal F}_{\mu-\nu_1}+\ldots+\hat{g}_k{\cal F}_{\mu-\nu_k})&\to 0\\
\end{array}$$
\normalsize

For $s>0$ and every $\mu'$ we have $H^s(D,{\cal F}_{\mu'})\simeq H^s(X^{alg},{\cal F}^{alg}_{\mu'})$:\\
Note that $R^s\pi_*{\cal F}_{\mu'}$ is concentrated upon $0$, so\\
$H^s(D,{\cal F}_{\mu'})\simeq (R^s\pi_*{\cal F}_{\mu'})_0\simeq (R^s\pi_*^{alg}{\cal F}^{alg}_{\mu'})_0\simeq H^0(\mathbb{C}^n,R^s\pi_*^{alg}{\cal F}^{alg}_{\mu'})\simeq H^s(X^{alg},{\cal F}^{alg}_{\mu'})$.\\
Now we can calculate $H^s(X^{alg},{\cal F}^{alg}_{\mu'})$ using [KKMS] or [O], as we will see.\\
If $h_{\mu'}$ is convex we have $H^s(X^{alg},{\cal F}^{alg}_{\mu'})=0$:\\
This is clear by [KKMS] Cor. 2, p. 44 or [O] Theorem 2.7 if $h_{\mu'}(\mathbb{N}^n)\subset\mathbb{Z}$. Otherwise there is an $M\in\mathbb{N}\setminus\{0\}$ such that $h_{M\mu'}=Mh_{\mu'}$ has this property. Now look at the Galois covering $(\mathbb{C}^*)^n\to (\mathbb{C}^*)^n: (z_1,\ldots,z_n)\mapsto (z_1^M,\ldots,z_n^M)$ with Galois group $G:=\{c\in\mathbb{C}^n\,|\,c^M={\bf 1}\}\simeq(\mathbb{Z}/M\mathbb{Z})^n$. The mapping extends to a Galois covering $p:X\to X$ with the same Galois group $G$. Now $G$ acts on the sheaf $p^{alg}_*{\cal F}^{alg}_{M\mu'}$ which admits a corresponding decomposition $p^{alg}_*{\cal F}^{alg}_{M\mu'}=\oplus_{c\in G}\,p^{alg}_*{\cal F}^{alg}_{M\mu'}(c)$. For $c={\bf 1}$ we obtain the $G$-invariant part: $p^{alg}_*{\cal F}^{alg}_{M\mu'}({\bf 1})=(p^{alg}_*{\cal F}^{alg}_{M\mu'})^{inv}\simeq {\cal F}^{alg}_{\mu'}$. We obtain a corresponding action on the cohomology groups, so: $(H^s(X^{alg},p^{alg}_*{\cal F}^{alg}_{M\mu'}))(c)\simeq H^s(X^{alg},(p^{alg}_*{\cal F}^{alg}_{M\mu'})(c))$. For $c={\bf 1}$ we obtain: $(H^s(X^{alg},p^{alg}_*{\cal F}^{alg}_{M\mu'}))^{inv}\simeq H^s(X^{alg},(p^{alg}_*{\cal F}^{alg}_{M\mu'})^{inv})\simeq H^s(X^{alg},{\cal F}^{alg}_{\mu'})$. Finally, $H^s(X^{alg},p^{alg}_*{\cal F}^{alg}_{M\mu'})=H^s(X^{alg},{\cal F}^{alg}_{M\mu'})=0$, hence $H^s(X^{alg},{\cal F}^{alg}_{\mu'})=0$ for $s>0$.\\ 
Furthermore, $H^0(D,{\cal F}_{\mu'})\simeq (\pi_*{\cal F}_{\mu'})_0\simeq ((\pi^{alg}_*{\cal F}^{alg}_{\mu'}g)_0)^{an}\simeq H^0(\mathbb{C}^n,\pi^{alg}_*{\cal F}^{alg}_{\mu'})\otimes_{\mathbb{C}[z_1,\ldots,z_n]}{\cal O}_{\mathbb{C}^n,0}$, and $H^0(\mathbb{C}^n,\pi^{alg}_*{\cal F}^{alg}_{\mu'})\simeq H^0(X^{alg},{\cal F}^{alg}_{\mu'})$. The latter is the ideal in $\mathbb{C}[z_1,\ldots,z_n]$ generated by all $z^q$ with $<p_j,q>\ge\mu'_j$ for all $j$. Altogether we obtain:\\
$H^0(D,{\cal F}_{\mu'})=F_{\mu'}({\cal O}_{\mathbb{C}^n,0})$.\\ 
We apply this to $\mu':=\mu-\nu_{j_1}-\ldots-\nu_{j_s}$, $\mu$ as in the hypothesis, $1\le j_1<\ldots<j_s\le k$. Note that $h_{\mu'}$ is convex, because $h_{\mu'}=h_{\mu-\nu_1-\ldots-\nu_k}+h_{\nu_{j_{s+1}}+\ldots+\nu_{j_k}}$ with $\{1,\ldots,k\}=\{j_1,\ldots,j_k\}$, and the latter two functions are convex.\\
So the vertical arrows are isomorphisms except possibly for the last one.\\
Using Lemma 2.2 we conclude that the lower row is exact: \\
Let us look at the complex of sheaves ${\cal G}^\cdot$:
$$\ldots 0\to {\cal F}_{\mu-\nu_1-\ldots-\nu_k}|D\to \oplus {\cal F}_{\mu-\nu_1-\ldots-\hat{\nu}_j-\ldots-\nu_k}|D\to\ldots\to\oplus {\cal F}_{\mu-\nu_j}|D\to {\cal F}_\mu|D\to 0\ldots$$
which constitutes a left resolution of $({\cal F}_{\mu}|D)/(\hat{g}_1{\cal F}_{\mu-\nu_1}|D+\ldots+\hat{g}_k{\cal F}_{\mu-\nu_k}|D$).\\
So ${\cal H}^l({\cal G}^\cdot)=0$ for $l\neq k$, ${\cal H}^k({\cal G}^\cdot)=({\cal F}_{\mu}|D)/(\hat{g}_1{\cal F}_{\mu-\nu_1}|D+\ldots+\hat{g}_k{\cal F}_{\mu-\nu_k}|D)$.\\
Therefore $\mathbb{H}^l(D,{\cal G}^\cdot)=0, l<k$, $\mathbb{H}^k(D,{\cal G}^\cdot)=H^0(D,{\cal F}_{\mu}/(\hat{g}_1{\cal F}_{\mu-\nu_1}+\ldots+\hat{g}_k{\cal F}_{\mu-\nu_k}))$.\\
On the other hand, we have seen that the sheaves ${\cal G}^l$ are acyclic, so:\\
$\mathbb{H}^l(D,{\cal G}^\cdot)=H^l(H^0(D,{\cal G}^\cdot))$ for all $l$.\\
Altogether we obtain the exactness of the lower row.\\
So the upper row is exact up to $F_\mu$. It is clear then that the whole upper row is exact, so the last vertical arrow has to be an isomorphism, too. This shows (*).

\vspace{2mm}
Now we show the same with $\mu+{\bf 1}$ instead of $\mu$.\\
Choose $\mu'$ such that $h_{\mu'}$ is convex.\\
If $\sigma$ is a cone of the fan as above, ${\cal F}_{\mu'+{\bf 1}}(X_\sigma)$ is generated by all $z^q$ with $q_{i_1}\ge 0,\ldots,q_{i_l}\ge 0, <p_{j_1},q>>\mu'_{j_1},\ldots,<p_{j_s},q>>\mu'_{j_s}$.\\
This means: $<p,q>\ge h_{\mu'}(p)+\eta h_{\nu_1}(p)$ if $\eta\in\mathbb{Q}$ is chosen so that $0<\eta\ll 1$. Note that $h_{\nu_1}$ is convex, so $h_{\mu'}+\eta h_{\nu_1}$, too, so we get $H^s(X^{alg},{\cal F}_{\mu'+{\bf 1}})=0, s>0$, as above.
 
\vspace{2mm}
Finally (**) follows from Lemma 2.3a); note that we have shown above the bijectivity of the last vertical arrow.

\vspace{2mm}
b) follows from a). Note that it was shown above that $F_\mu/(g_1F_{\mu-\nu_1}+\ldots+g_kF_{\mu-\nu_k})\simeq H^0(D,{\cal F}_\mu/\hat{g}_1{\cal F}_{\mu-\nu_1}+\ldots+\hat{g}_k{\cal F}_{\mu-\nu_k})$; cf. Lemma 2.3b).\\

{\bf 3. One-index Newton filtration}\\

The methods introduced at the beginning can also be used in connection with the results of Kushnirenko [K] as well as Bivi\`a-Ausina, Fukui and Saia [BFS] about the one-index Newton filtration.\\
Here we need a more serious assumption about the Newton polyhedra than before. For instance, we may suppose that $\nu_1,\ldots,\nu_k$ are scalar multiples of each other, so $\nu_i=\rho_i\nu$ for some $\rho_i>0, \nu\in\mathbb{N}^n$. Then put $\Delta:=\{q\,|\,<p,q>\ge\nu_j,j=1,\ldots,r\}$.\\
We gain some generality, however, if we follow the proposal made in [BFS] and start from $\Delta$.\\
So let $\Delta$ be a Newton polyhedron in the following sense: $\Delta$ is the convex hull of $S\cup\mathbb{R}_+^n$, where $S$ is some subset of $\mathbb{N}^n$ which intersects each of the coordinate axes. The compact facets of $\Delta$ are again of the form $\{<p_j,q>\ge \nu_j\}$, $j=1,\ldots,r$, with $p_j\in\mathbb{N}^n$. Let us fix a simplicial fan with edges generated by $e_1,\ldots,e_n,p_1,\ldots,p_r$, and let $X$ be the corresponding toric variety.\\
Choose $M\in\mathbb{N}\setminus\{0\}$ such that $\psi:\mathbb{N}^n\to\mathbb{R}$: $\psi(q):=M\,\min\{\frac{1}{\nu_j}<p_j,q>\,|\,j=1,\ldots,r\}$ takes values in $\mathbb{Z}$.\\
Then let $\hat{F}_l$ be the ideal in ${\cal O}_{\mathbb{C}^n,0}$ generated by all $z^q$ such that $\psi(q)\ge l$, i.e. $<p_j,q>\ge\frac{l}{M}\nu_j$ for all $j$. (We could take $M$ minimal but this is not important).  We can introduce corresponding sheaves $\hat{\cal F}_l$ on $X$, similarly as in \S 2.\\
Choose $h_l:\mathbb{R}_+^n\to\mathbb{R}$ such that $h_l(e_i)=0, i=1,\ldots,n$, $h_l(p_j)=\frac{l\nu_j}{M}$, $j=1,\ldots,r$, $h_l|\sigma$ linear for all cones of the fan. Then $h_l$ is convex, $l\ge 0$, because $h_{Ml}$ is convex, $h_l(\mathbb{N}^n)\subset \mathbb{Z}$ if $M|l$.\\
Now fix $\rho_i\in\mathbb{N}\setminus\{0\}$ and $g_i\in \hat{F}_{\rho_i}$, $i=1,\ldots,k$.\\
If $\sigma$ is a cone of the fan, generated by $e_{i_1},\ldots,e_{i_l},p_{j_1},\ldots,p_{j_s}$, let $in_\sigma g_i$ be the part of $g_i$ involving only those $z^q$ with $q_{i_1}\ge 0,\ldots,q_{i_l}\ge 0, <p_{j_1},q>\ge \frac{\nu_{j_1}\rho_i}{M},\ldots,<p_{j_s},q>\ge\frac{\nu_{j_s}\rho_i}{M}$.\\

{\bf Lemma 3.1:} Suppose that we have the same hypothesis as in Lemma 2.2.\\
a) There is a commutative diagram with exact rows
$$\begin{array}{cccc}
0\to &\hat{F}_{l-\rho_1-\ldots-\rho_k}&\to\ldots\to&\hat{F}_l\\
&\downarrow\simeq&&\downarrow\simeq\\
0\to &H^0(D,\hat{\cal F}_{l-\rho_1-\ldots-\rho_k})&\to\ldots\to&H^0(D,\hat{\cal F}_l)
\end{array}$$
b) There is an exact sequence 
$$0\to \hat{F}_{l-\rho_1-\ldots-\rho_k}\to\ldots\to \hat{F}_l\to \hat{F}_l/(g_1\hat{F}_{l-\rho_1}+\ldots+g_k\hat{F}_{l-\rho_k})\to 0$$

{\bf Proof:} It is sufficient to show a).\\
We proceed as in the proof of Theorem 2.7. In fact we will show that $H^i(D,\hat{\cal F}_l)=0$, $i\neq 0,n-1$: \\
We concentrate on the case $i\neq 0$. As in the proof of Theorem 2.7, we have $H^i(D,\hat{\cal F}_l)\simeq H^i(X^{alg},\hat{\cal F}^{alg}_l)$. We will show that again it is sufficient to look at the case $M|l$.\\
Look at the Galois covering $(\mathbb{C}^*)^n\to(\mathbb{C}^*)^n$: $(z_1,\ldots,z_n)\mapsto (z_1^M,\ldots,z_n^M)$ with Galois group $G:=\{c\in\mathbb{C}^n\,|\,c^M={\bf 1}\}\simeq (\mathbb{Z}/(M\mathbb{Z}))^n$. The mapping extends to a Galois covering $p:X\to X$ (with Galois group $G$, too, of course). As in the proof of Theorem 2.7, we have an action on $p^{alg}_*\hat{\cal F}^{alg}_{Ml}$ and its cohomology groups, we have $(p^{alg}_*\hat{\cal F}^{alg}_{Ml})^{inv}\simeq \hat{\cal F}^{alg}_l$, and $H^i(X^{alg},\hat{\cal F}^{alg}_l)\simeq H^i(X^{alg},(p^{alg}_*\hat{\cal F}^{alg}_{Ml})^{inv})\simeq (H^i(X^{alg},\hat{\cal F}^{alg}_{Ml}))^{inv}$, so it is sufficient to look at the case $M|l$.\\
So we suppose from now on that $M|l$. The advantage is that $\hat{\cal F}_l$ is then a line bundle.\\
Now if $l\ge 0$ we know that $h_l$ is convex, so $H^i(X^{alg},\hat{\cal F}^{alg}_l)=0$.\\
So assume $l<0$. Then $lh$ is concave, with $h:=h_1$. This implies for $q\in\mathbb{Z}^n$ the following: the set $A^l_q:=\{p\in |X|\,|\,<p,q>\ge lh(p)\}$ is convex.\\
Let $i>0$. As in the proof of Theorem 2.7, $H^i(D,\hat{\cal F}_l)=H^i(X^{alg},\hat{\cal F}^{alg}_l)$. By [KKMS] p. 42, we have for the part corresponding to the character $q$:\\
$H^i(X^{alg},\hat{\cal F}^{alg}_l)(q)=H^i(|X|,|X|\setminus A_q^l;\mathbb{C})$. Here $|X|=\mathbb{R}^n_+$.\\
Take the affine hypersurface $H=\{t_1+\ldots+t_n=1\}$ in $\mathbb{R}^n$. Then:\\
$H^i(|X|,|X|\setminus A_q^l;\mathbb{C})=H^i(|X|\cap H,|X|\cap H\setminus A_q^l;\mathbb{C})$,\\
because $A_q^l$ is conic, i.e. invariant under $p\mapsto \lambda p$, $\lambda\ge 0$ (note that $h(\lambda p)=\lambda h(p)$). So $|X|\cap H\setminus A_q^l$ is a deformation retract of $|X|\setminus A^l_q$.\\
As we have seen, $A^l_q$ is convex, so $A^l_q\cap H$, too.\\
Now $h|(\partial|X|)\equiv 0$ because the cones contained in $\partial|X|$ are generated by standard unit vectors. If such a cone is generated by $e_{i_1},\ldots,e_{i_l}$, we have that $p\in\sigma$ belongs to $A^l_q$ if and only if $p_{i_1}q_{i_1}+\ldots+p_{i_l}q_{i_l}\ge 0$. Therefore $A^l_q\cap H\cap\partial|X|$ is contractible or empty as soon as $A^l_q\cap H\cap\partial|X|$ is a proper subset of $H\cap\partial|X|$.\\
Let $C^l_q$ be the set of all $p\in H\setminus|X|$ which lie on a line joining the point $(\frac{1}{n},\ldots,\frac{1}{n})$ with some point of $A^l_q\cap H\cap\partial|X|$. Consider cohomology with complex coefficients. Then we have: $H^i(|X|\cap H,|X|\cap H\setminus A^l_q)\simeq H^i(H,H\setminus A^l_q\cup C^l_q)$. Look at the long exact sequence:\\
$\ldots\to H^i(H,H\setminus A^l_q)\to H^i(H,H\setminus A^l_q\cup C^l_q)\to H^i(H\setminus A^l_q,H\setminus A^l_q\cup C^l_q)\to\ldots$\\
We have $H^i(H,H\setminus A^l_q)\simeq H_{n-1-i}(A^l_q\cap H)$, so $H^i(H,H\setminus A^l_q) =0$ for $i\neq n-1$. Moreover,
$H^{n-1}(H,H\setminus A^l_q)\simeq\mathbb{C}$ if $A^l_q\neq \{0\}$ and $= 0$ otherwise. Note that $0\in A^l_q$.\\
On the other hand, $H^i(H\setminus A^l_q,H\setminus A^l_q\cup C^l_q)\simeq H^i(H\setminus |X|,H\setminus |X|\cup C^l_q)\simeq H^i(H\cap\partial|X|,H\cap\partial|X|\setminus A^l_q)\simeq H_{n-2-i}(A^l_q\cap\partial|X|\cap H)$. Assume that $A_l^q\cap H\cap\partial|X|$ is a proper subset of $H\cap\partial|X|$, i.e. $A^l_q\neq |X|$. Then $A_l^q\cap H\cap\partial|X|$ is contractible or empty, as we saw above. So $H^i(H\setminus A^l_q,H\setminus A^l_q\cup C^l_q)=0$ for $i\neq n-2$, $H^{n-2}(H\setminus A^l_q,H\setminus A^l_q\cup C^l_q)\simeq\mathbb{C}$ if $A^l_q\cap \partial|X|\neq\{0\}$ and $=0$ otherwise.\\
Finally, note that the mapping $H^i(H\setminus A^l_q,H\setminus A^l_q\cup C^l_q)\to H^{i+1}(H,H\setminus A^l_q)$ is bijective for all $i$ if $A^l_q\neq |X|, A^l_q\cap\partial|X|\neq\{0\}$:\hfill(*)\\
Choose $p\in A^l_q\cap H\cap\partial|X|$ and let $C^*_p$ be the set of points of the line through $(\frac{1}{n},\ldots,\frac{1}{n})$ and $p$ which lie outside the interior of $|X|$.  Now look at the commutative diagram
$$\begin{array}{ccc}
H^i(H\setminus\{p\},H\setminus C^*_p)&\to& H^{i+1}(H,H\setminus \{p\})\\
\downarrow&&\downarrow\\
H^i(H\setminus A^l_q,H\setminus A^l_q\cup C^l_q)&\to& H^{i+1}(H,H\setminus A^l_q)\\
\end{array}$$
The right vertical arrow is bijective: this follows from the commutative diagram
$$\begin{array}{ccc}
H^{i+1}(H,H\setminus \{p\})&\to&H_{n-1-i}(\{p\})\\
\downarrow&&\downarrow\\
H^{i+1}(H,H\setminus A^l_q)&\to&H_{n-1-i}(A^l_q\cap H)
\end{array}$$
where the horizontals and the right vertical are bijective.\\
Similarly, the left vertical of the first diagram is bijective, too, as well as the upper horizontal: we have that $H^i(H,H\setminus C_p^*)=0$ for all $i$. So (*) is proved.\\
Collecting these results we see the following. We distinguish four cases:\\
a) $A^l_q=\{0\}$: then $H^i(|X|,|X|\setminus A^l_q)=0$ for all $i$.\\
b) $A^l_q\neq\{0\}$, $A^l_q\cap\partial|X|=\{0\}$: Then $H^i(|X|,|X|\setminus A^l_q)=0$ for $i\neq n-1$, $H^{n-1}(|X|,|X|\setminus A^l_q)\simeq\mathbb{C}$.\\
c) $A^l_q\neq|X|, A^l_q\cap\partial|X|\neq\{0\}$: Then $H^i(|X|,|X|\setminus A^l_q)=0$ for all $i$.\\
d) $A^l_q=|X|$: Then $H^i(|X|,|X|\setminus A^l_q)=0$ for $i\neq 0$, $H^0(|X|,|X|\setminus A^l_q)\simeq\mathbb{C}$.\\
Altogether, we obtain that $H^i(|X|,|X|\setminus A^l_q)=0$ for $i\neq 0,n-1$, so $H^i(D,\hat{\cal F}_l)=0$ for $i\neq 0,n-1$.

\vspace{2mm}
Let ${\cal G}^\cdot$ be similarly defined as in the proof of Theorem 2.7. As shown loc.cit. we have: $\mathbb{H}^i(D,{\cal G}^\cdot)=0$, $i<k$. We have a spectral sequence which converges to this hypercohomology with $E_1^{pq}=H^q(D,{\cal G}^p)$. By the cohomology calculation above we have $E_1^{pq}=0$ unless $p\in \{0,\ldots,k\}, q\in\{0,n-1\}$. Therefore we have a long exact sequence
$$\ldots\to \mathbb{H}^i(D,{\cal G}^\cdot)\to E_2^{i-n+1,n-1}\to E_2^{i+1,0}\to \mathbb{H}^{i+1}(D,{\cal G}^\cdot)\to\ldots$$
see [CE] XV Theorem 5.11. Now for $i<k\le n$ we have $E_2^{i-n,n-1}=0$, $\mathbb{H}^i(D,{\cal G}^\cdot)=0$, so $0=E_2^{i,0}=H^i(H^0(D,{\cal G}^\cdot))$. So we get that the lower row of the diagram in a) is exact:
$$\begin{array}{cccc}
0\to &\hat{F}_{l-\rho_1-\ldots-\rho_k}&\to\ldots\to&\hat{F}_l\\
&\downarrow&&\downarrow\\
0\to &H^0(D,\hat{\cal F}_{l-\rho_1-\ldots-\rho_k})&\to\ldots\to&H^0(D,\hat{\cal F}_l)
\end{array}$$
The vertical arrows are still isomorphisms.\\
Obviously this implies our lemma.\\

{\bf Remark:} Assume furthermore that $k\le n-2$. Then the last vertical arrow of the commutative diagram with exact upper row:
\footnotesize
$$\begin{array}{ccccccc}
0\to&\hat{F}_{l-\rho_1-\ldots-\rho_l}&\to\ldots\to&\hat{F}_l&\to&\hat{F}_l/(g_1\hat{F}_{l-\rho_1}+\ldots+g_k\hat{F}_{l-\rho_k})&\to 0\\
&\simeq\downarrow&&\simeq\downarrow&&\downarrow&\\
0\to&H^0(D,\hat{\cal F}_{l-\rho_1-\ldots-\rho_l})&\to\ldots\to&H^0(D,\hat{\cal F}_l)&\to&H^0(D,\hat{\cal F}_l/(\hat{g}_1\hat{\cal F}_{l-\rho_1}+\ldots+\hat{g}_k\hat{\cal F}_{l-\rho_k}))&\to 0\\
\end{array}$$
\normalsize
is bijective, i.e. the lower row is exact, too.\\
b) If $k=n-1$, the last vertical of the commutative diagram with exact upper row:
\footnotesize
$$\begin{array}{ccccccc}
0\to&\hat{F}_{l-\rho_1-\ldots-\rho_l}&\to\ldots\to&\hat{F}_l&\to&\hat{F}_l/(g_1\hat{F}_{l-\rho_1}+\ldots+g_k\hat{F}_{l-\rho_k})&\to 0\\
&\simeq\downarrow&&\simeq\downarrow&&\downarrow&\\
0\to&H^0(D,\hat{\cal F}_{l-\rho_1-\ldots-\rho_l})&\to\ldots\to&H^0(D,\hat{\cal F}_l)&\to&H^0(D,\hat{\cal F}_l/(\hat{g}_1\hat{\cal F}_{l-\rho_1}+\ldots+\hat{g}_k\hat{\cal F}_{l-\rho_k}))&\\
\end{array}$$
\normalsize
is injective, i.e. the lower row is exact, too.\\

{\bf Proof:} We start from the proof of Lemma 3.1. We have furthermore: $\mathbb{H}^k(D,{\cal G}^\cdot)=H^0(D,\hat{\cal F}_l/(\hat{g}_1\hat{\cal F}_{l-\rho_1}+\ldots+\hat{g}_k\hat{\cal F}_{l-\rho_k}))$, see proof of Theorem 2.7.\\
On the other hand, $E_2^{k0}=H^k(H^0(D,{\cal G}^\cdot))=coker(H^0(D,\oplus\hat{\cal F}_{\l-\rho_j})\to H^0(D,\hat{\cal F}_l))\simeq coker(\oplus\hat{F}_{l-\rho_j}\to\hat{F}_l)=\hat{F}_l/(g_1\hat{F}_{l-\rho_1}+\ldots+g_k\hat{F}_{l-\rho_k})$. Finally, we have an exact sequence
$$E_2^{k-n,n-1}\to E_2^{k,0}\to \mathbb{H}^k(D,{\cal G}^\cdot)\to E_2^{k-n+1,n-1}$$
For $k\le n-2$ we have $E_2^{k-n,n-1}=E_2^{k-n+1,n-1}=0$, for $k=n-1$ we have $E_2^{k-n,n-1}=0$.\\ 

{\bf Lemma 3.2:} Suppose that we have the same hypothesis as in Lemma 2.3.\\
a) We have an exact sequence  
$0\to \hat{F}_{l-\rho_1-\ldots-\rho_k}/\hat{F}_{l-\rho_1-\ldots-\rho_k+1}\to\ldots\to \hat{F}_l/\hat{F}_{l+1}\to\hat{F}_l/(\hat{F}_{l+1}+g_1\hat{F}_{l-\rho_1}+\ldots+g_k\hat{F}_{l-\rho_k})\to 0$.\\
b) We have $\hat{F}_l/(g_1\hat{F}_{l-\rho_1}+\ldots+g_k\hat{F}_{l-\rho_k})\simeq \hat{F}_l/((g_1,\ldots,g_k)\cap \hat{F}_l)$.\\

{\bf Proof:} a) Suppose first that $k\le n-1$. We have a commutative diagram
$$\begin{array}{ccc} \hat{F}_{l+1}/(g_1\hat{F}_{l+1-\rho_1}+\ldots+g_k\hat{F}_{l+1-\rho_k})&\to&\hat{F}_l/(g_1\hat{F}_{l-\rho_1}+\ldots+g_k\hat{F}_{l-\rho_k})\\
\downarrow&&\downarrow\\
H^0(D,\hat{\cal F}_{l+1}/(\hat{g}_1\hat{F}_{l+1-\rho_1}+\ldots+\hat{g}_k\hat{F}_{l+1-\rho_k}))&\to&H^0(D,\hat{\cal F}_l/(\hat{g}_1\hat{F}_{l-\rho_1}+\ldots+\hat{g}_k\hat{F}_{l-\rho_k}))
\end{array}$$
The lower row is injective by Lemma 2.3. The vertical arrows are injective by the preceding Remark. Therefore the upper row must be injective, too. So we obtain our result by a short exact sequence of sheaf complexes.\\
The case $k=n$ is covered by [BFS] p. 72; in the special case treated by Kushnirenko see [K] Th\'eor\`eme 2.8.\\
b) We must show that the sequence $\oplus \hat{F}_{l-\rho_j}\to\hat{F}_l\to\hat{F}_l/(g_1,\ldots,g_k)\cap\hat{F}_l$ is exact. We proceed similarly as in the proof of Lemma 1.2. Assume that $h\in\hat{F}_l$ represents an element of the kernel of the second map. Then we can write $h=\sum g_ih_i$ with $h_i\in{\cal O}_{\mathbb{C}^n,0}$. Choose $l'\le l$ maximal such that $h_i\in\hat{F}_{l'-\rho_i}$ for all $i$. We must show that $l'=l$. Suppose that $l'<l$: By a) we have an exact sequence
$\oplus_{i<j}\hat{F}_{l'-\rho_i-\rho_j}/\hat{F}_{l'+1-\rho_i-\rho_j}\to \oplus_i\hat{F}_{l'-\rho_i}/\hat{F}_{l'+1-\rho_i}\to\hat{F}_{l'}/\hat{F}_{l'+1}$.
Now $(h_1,\ldots,h_k)$ represents an element of the kernel of the second map, so we can find an inverse image in the first group which can be used in order to modify $(h_1,\ldots,h_k)$ in such a way that we have $l'+1$ instead of $l'$, contradiction.\\

{\bf Remark:} For $k\le n-1$, $l\ge\rho_1+\ldots+\rho_k$ the proof can be simplified, as in the proof of Theorem 2.7. \\

{\bf Lemma 3.3:} Let $P_l$ be the coefficient of $P_{{\cal O}_\mathbb{C}^n,0}(t)$ at $t_1^l\cdot\ldots\cdot t_r^l$. Then $\hat{P}(\tau):=\sum_lP_l\tau^l$ is the Poincar\'e series of ${\cal O}_{\mathbb{C}^n,0}$ with respect to the one-index filtration $\hat{F}$.\\
Furthermore, let $M_l$ be the set of all $q\in\mathbb{N}^n$ such that $<p_j,q>\ge l$ for all $j$. Then:\\
$\hat{P}(\tau)=\sum_l\#(M_l\setminus M_{l+1})\tau^l$.\\

The proof is straightforward. \\

{\bf Theorem 3.4:} Put $\hat{Q}(\tau):=(1-\tau^{\rho_1})\cdot\ldots\cdot(1-\tau^{\rho_k})\hat{P}(\tau)=\sum_l\hat{Q}_l\tau^l$.\\
a) Suppose that the hypothesis of Lemma 2.3 is fulfilled. Then the Poincar\'e series of ${\cal O}_{Y,0}$ with respect to the one-index filtration $\hat{F}_l$ is $\hat{Q}(\tau)$.\\
b) If $k=n$, $\dim\,{\cal O}_{Y,0}$ is the sum of coefficients of this series $\hat{Q}(\tau)$. So the latter is a polynomial, and $\dim\,{\cal O}_{Y,0}=\hat{Q}(1)$.\\

{\bf Proof:} a) follows from Lemma 3.2.\\
b) clear.\\

{\bf Remark:} a) We can also compute the Poincar\'e series with respect to the filtration $H^0(D,\hat{\cal F}_l/(\hat{g}_1,\ldots,\hat{g}_k)\cap\hat{\cal F}_l)$: For $k\le n-2$ this filtration coincides with $\hat{F}_l$, see Remark a) after Lemma 3.1. The case $k=n$ is uninteresting because of Lemma 2.6c). For $k=n-1$ the Poincar\'e series for $H^0(D,\hat{\cal F}_l/(\hat{g}_1,\ldots,\hat{g}_k)\cap\hat{\cal F}_l)$ may be computed with Theorem 2.4, since the upper cohomology groups vanish: the support of the sheaf is concentrated upon a one-dimensional set, its intersection with $D$ is zero-dimensional.\\
b) By Theorem 2.7, we get in the case $k=n$ that $\hat{Q}_l=0$ for $l\ge \rho_1+\ldots+\rho_n$; note that the third series there is $0$. Cf. Theorem 3.4: we know a priori at least that $\hat{Q}_l=0,l\gg 0$, $\hat{Q}$ is a polynomial.\\

{\bf 4. Appendix: Remarks on irreducibility}\\

Here we want to prove an assertion of [EG2] (see Remark after Lemma 2.6) and discuss the difficulties which arise when trying to generalize it.\\

Let $g$ be a non-degenerate Laurent polynomial in $n$ variables. Let $\Delta$ be the convex hull of $supp\,g$. Assume that $\Delta$ is full, i.e. $\dim\Delta=n$. \\

{\bf Theorem 4.1:} $H_j((\mathbb{C}^*)^n,(\mathbb{C}^*)^n\cap \{g=0\};\mathbb{Z})=0$ if $n\ge 2, j\le 1$.\\

{\bf Proof:} This follows from a Lefschetz theorem of Oka, see [Ok] V Cor. 4.6.1.\\

Let ${\cal F}_0$ be the dual fan to $\Delta$, $\cal F$ a refinement, and $X$ the corresponding toric variety. Then $\Delta$ corresponds to a line bundle on $X$, and $g$ can be regarded as a global section $\rm g$ in the latter. Let $Y$ be the zero locus of $\rm g$.\\

{\bf Corollary 4.2:} For $n\ge 2$, $Y$ is non-empty and irreducible, and $Y$ is the closure of $(\mathbb{C}^*)^n\cap\{g=0\}$ in $X$.\\

{\bf Proof:} Assume first that $X$ is smooth, so $Y$ is a smooth hypersurface in $X$. Theorem 4.1 implies that $H_0((\mathbb{C}^*)^n\cap Y;\mathbb{Z})\simeq\mathbb{Z}$, so $Y\neq\emptyset$, and $(\mathbb{C}^*)^n\cap Y$ is connected. Assume that $Y$ is not irreducible, i.e. not connected: then there is a connected component $Y_i$ which is contained in the complement of $(\mathbb{C}^*)^n$. Furthermore, $\dim\,Y_i=n-1$, so there is a $(\mathbb{C}^*)^n$-orbit $O$ in $X\setminus(\mathbb{C}^*)^n$ such that $O\cap Y_i\neq \emptyset$ and $\dim\,O\cap Y_i=n-1$. But this is impossible because $Y_i$ is a hypersurface in $O$ and $\dim\,O\le n-1$, hence $\dim\,O\cap Y_i\le n-2$.\\
If $X$ is not smooth there is a modification $\pi:\tilde{X}\to X$, $\tilde{X}$ smooth toric. Let $\tilde{Y}$ be correspondingly defined, then $\tilde{Y}$ is irreducible, hence $Y=\pi(\tilde{Y})$, too.\\
The rest is clear.\\

{\bf Remark:} a) The assumption that $\Delta$ is full is unnecessary if we replace the condition $n\ge 2$ by $\dim\,\Delta\ge 2$. Use [Ok] V Th. 4.6 here. Note that ${\cal F}_0$ is no longer a fan but at least a cone decomposition.\\
b) For the statement about $Y$ being the closure the assumption $n\ge 2$ is not necessary.\\

Now let $g$ be a convenient polynomial in $n$ variables which is non-degenerate at $0$. Let $\Delta$ be the convex hull of $supp\,g+\mathbb{N}^n$. Let $\cal F$ be a refinement of the dual fan to $\Delta$ and $X$ the corresponding toric variety. Then we have a toric modification $\pi:X\to \mathbb{C}^n$. Let $Y$ be the strict transform of $g^{-1}(\{0\})$.\\

{\bf Theorem 4.3:} (cf. [EG2]) If $D$ is a component of the exceptional divisor which corresponds to a compact facet of $\Delta$, $n\ge 3$, $D\cap Y$ is an irreducible hypersurface in $D$.\\

{\bf Proof:} Let $\sigma$ be the corresponding facet. Then $\sigma$ is full with respect to the hyperplane which contains $\sigma$. So $\sigma$ corresponds to an $(n-1)$-dimensional toric variety $D$ - the closure of the orbit which corresponds to $\sigma$. Note that $\dim\,D\ge 2$. Let $\tau$ be the edge which is dual to $\sigma$. Then $in_\tau g$ defines a smooth non-degenerate hypersurface in $D$, namely $Y\cap D$. So we can apply Theorem 4.1.\\

The attempt to pass from hypersurfaces to complete intersections involves heavy restrictions.\\

Let $g_1,\ldots,g_k$ be Laurent polynomials which define a non-degenerate complete intersection in $(\mathbb{C}^*)^n$, see [Ok] p. 107. Let $\Delta_i$ be the Newton polyhedron of $g_i$. Then we have the following generalization of Theorem 4.1:\\

{\bf Theorem 4.4}: Assume that $\Delta_1,\ldots,\Delta_k$ are full. Then $H_j((\mathbb{C}^*)^n,(\mathbb{C}^*)^n\cap\{g_1=\ldots=g_k=0\};\mathbb{Z})=0$ if $n-k\ge 1, j\le 1$.\\

{\bf Proof:} Use [Ok] V Cor. 4.6.1 again.\\

Let $\cal F$ be a fan which refines the dual fan to $\Delta:=\Delta_1+\ldots+\Delta_k$ and $X$ the corresponding toric variety. Note that $\Delta_i$ correesponds to a line bundle on $X$, and we have a section ${\rm g}_i$ in it which corresponds to $g_i$. Then ${\rm g}_1,\ldots,{\rm g}_k$ define a complete intersection $Y$ in $X$. \\

{\bf Corollary 4.5:} For $n-k\ge 1$, $Y$ is non-empty and irreducible, and $Y$ coincides with the closure of $(\mathbb{C}^*)^n\cap\{g_1=\ldots=g_k=0\}$ in $X$.\\

The proof is analogous to the one of Corollary 4.2.\\

{\bf Remark:} We may weaken the assumption that all Newton polyhedra are full: Assume that $\dim(\Delta_1+\ldots+\Delta_j)=\dim\;\Delta_j$, $2\le j\le k$, and 
$\dim\Delta_j>j, j=1,\ldots,k$. Then apply [Ok] V Theorem 4.6.\\
We get a corresponding modification of Corollary 4.5.\\

Now let us pass to a generalization of Theorem 4.3.\\

{\bf Theorem 4.6:} Assume that the set of edges (i.e. one-dimensional cones) of the dual fan to $\Delta_i$ does not depend on $i$. (This holds, for example, if $\Delta_1,\ldots,\Delta_k$ are similar to each other, or if the dual fans to $\Delta_1,\ldots,\Delta_k$ coincide.) Then we get the assertion of Theorem 4.3 provided that we replace the assumption $n\ge 3$ by $n-k\ge 2$.\\

{\bf Proof:} Suppose that $\sigma$ is a corresponding facet. For each $i$ there is a facet $\sigma_i$ of $\Delta_i$ such that $\sigma$ and $\sigma_i$ are contained in parallel hypersurfaces. So we can pass to the toric variety $\overline{O}_\tau$, where $\tau$ is dual to $\sigma$, and study the complete intersection defined by $in_\tau g_1,\ldots,in_\tau g_k$. The corresponding Newton polyhedra $\sigma_1,\ldots,\sigma_k$ are full - so we can apply Corollary 4.5.\\

It is difficult to weaken the hypothesis of Theorem 4.6. Anyhow we have a partial result for the case $k=2$:\\

{\bf Remark:} Let $k=2$ and assume that $n\ge 4$ and that for every facet $\sigma_2$ of $\Delta_2$ of the form $\{q\in\Delta_2\,|\,<p,q>\ge s_2\}$ we have that $\dim\,\sigma_1\ge 2$, where $\sigma_1=\{q\in \Delta_1\,|\,<p.q>\ge s_1\}$ and $s_1=\min\{<p.q>\,|\,q\in\Delta_1\}$, and vice versa. Then we get the assertion of Theorem 4.6.\\

For the proof note the following: Let us fix facet $\sigma$ of $\Delta$. Without loss of generality we can assume that there is a facet $\sigma_2$ of $\Delta_2$ with the same normal vectors. Let $\sigma_1$ be as above and $\tau$ dual to $\sigma$. Then we pass to $\overline{O}_\tau$. The Newton polyhedra of $in_\tau g_1$ and $in_\tau g_2$ are $\sigma_1$ resp. $\sigma_2$, so the dimension is $\ge 2$ resp. $n-1\ge 3$. By the Remark after Corollary 4.5 we obtain our statement.\\

{\bf Example:} Put $g_1(z):=z_1+z_2+z_3+z_4^2+z_5^2+z_6^2$, $g_2(z):=z_1^2+z_2^2+z_3^2+z_4+z_5+z_6$. Then the assumption of the last Remark is fulfilled.

\vspace{2cm}

{\bf References}\\

[BFS] C. Bivi\`a-Ausina, T. Fukui, M.J. Saia: Newton filtrations, graded algebras and codimension of non-degenerate ideals. Math. Proc. Cambridge Philos. Soc. {\bf 133}, 55-75 (2002).\\

[CDG] A. Campillo, F. Delgado, S.M. Gusein-Zade: Multi-index filtrations and generalized Poincar\'e series. Monatsh. Math. {\bf 150}, 193-209 (2007).\\

[CE] H. Cartan, S. Eilenberg: Homological Algebra. Princeton Univ. Press: Princeton, N.J. 1956.\\

[EG1] W. Ebeling, S.M. Gusein-Zade: Multi-variable Poincar\'e series associated with Newton diagrams. J. of Sing. {\bf 1}, 60- 68 (2010).\\

[EG2] W. Ebeling, S.M. Gusein-Zade: On divisorial filtrations associated with Newton diagrams. J. of Sing. {\bf 3}, 1-7 (2011).\\

[H] R. Hartshorne: Algebraic geometry. Springer-Verlag: New York 1977.\\

[K] A.G. Kouchnirenko: Poly\`edres de Newton et nombres de Milnor. Inventiones Math. {\bf 32}, 1-31 (1976).\\

[KKMS] G. Kempf, F. Knudsen, D. Mumford, B. Saint-Donat: Toroidal Embeddings I. Springer LN {\bf 339} (1973).\\

[L] A. Lemahieu: Poincar\'e series of embedded filtrations. arXiv:0906.4184v1 [math.AG], 2009\\

[O] T. Oda: Convex Bodies and Algebraic Geometry. Springer-Verlag: Berlin 1988.\\

[Ok] M.Oka: Non-degenerate complete intersection singularity. Hermann: Paris 1997.\\

\end{document}